\documentclass{amsart}
\usepackage{amsmath, amssymb}

\newcommand{\al}{{\alpha}}
\newcommand{\be}{{\beta}} 

\newcommand{\bbC}{\mathbb{C}}
\newcommand{\bbF}{\mathbb{F}}
\newcommand{\bbK}{\mathbb{K}}
\newcommand{\bbL}{\mathbb{L}}
\newcommand{\bbQ}{\mathbb{Q}}
\newcommand{\bbR}{\mathbb{R}}
\newcommand{\bbZ}{\mathbb{Z}}
\newcommand{\bx}{{\bf x}}
\newcommand{\bX}{{\bf X}}
\newcommand{\by}{{\bf y}}
\newcommand{\bY}{{\bf Y}}

\newcommand{\cA}{{\mathcal{A}}}
\newcommand{\cH}{{\mathcal{H}}}
\newcommand{\cL}{{\mathcal{L}}}
\newcommand{\cM}{{\mathcal{M}}}
\newcommand{\cN}{{\mathcal{N}}}
\newcommand{\cO}{{\mathcal{O}}}
\newcommand{\cP}{{\mathcal{P}}}
\newcommand{\cS}{{\mathcal{S}}}
\newcommand{\cU}{{\mathcal{U}}}

\newcommand{\disc}{{\rm disc}}
\newcommand{\lcm}{{\rm lcm}}

\newtheorem{corollary}{Corollary}
\newtheorem{lemma}{Lemma}
\newtheorem{theorem}{Theorem}

\theoremstyle{definition}
\newtheorem{definition}{Definition}

\newtheorem{note}{Note}


\begin{document}

\title[Full modules and norm form equations]{Modules with many non-associates and norm form equations with many families of solutions}

\author{Paul M Voutier}
\address{London, UK \\paul.voutier@gmail.com}

\dedicatory{To Wolfgang M. Schmidt, with warmest wishes and deepest admiration on his 80th birthday}
\begin{abstract}
For every number field $\bbK$, with $[\bbK:\bbQ] \geq 3$, we show that the number of
non-associates of the same norm in a full module in $\bbK$ does not depend only on
$\bbK$, but can also depend on the module itself.

As a corollary, the same can be true for the number of families of solutions of
degenerate norm form equations. So the uniform bound obtained by Schmidt for
the number of solutions in the non-degenerate case does not hold always here.

For three-variable norm forms not arising from full modules, we do obtain a
Schmidt-type bound for the number of families of solutions that, together with
the above result, completes this aspect of the study of three-variable norm forms.
\end{abstract}

\maketitle 

\section{Introduction}

\subsection{Non-associates}
Let $\bbK$ be an algebraic number field with $r=[\bbK:\bbQ]$. Let $\al_{1}$,
\ldots, $\al_{n}$ lie in $\bbK$ and put $L(\bX) = \al_{1}X_{1} + \cdots + \al_{n}X_{n}$.

The set $\cM = \left\{ L(\bx): \bx \in \bbZ^{n} \right\}$ is a 
$\bbZ$-module contained in $\bbK$.

For every subfield $\bbL$ of $\bbK$, let $\cM^{\bbL}$ consist of 
the elements $\be$ of $\cM$ such that for every $\al \in \bbL$ 
there is a non-zero rational integer $z$ with $z \al \be \in \cM$. 

\begin{definition} 
We can associate with each $\cM^{\bbL}$ a {\it ring of coefficients}, 
which we will denote by $\cO_{\cM}^{\bbL}$, i.e., the set of 
$\al \in \bbL$ such that $\al \be \in \cM^{\bbL}$ for every 
$\be \in \cM^{\bbL}$.
\end{definition} 

For our purposes here, we single out a particular subgroup
of the group of units in $\bbL$: let $\cU_{\cM}^{\bbL}$ be the group of
elements in $\cO_{\cM}^{\bbL}$ of norm $1$.

\begin{definition} 
We say that $\cM$ is a {\it full module} in $\bbK$ if its rank, as a $\bbZ$-module,
is equal to $r$. 

Two elements $\mu_{1}$ and $\mu_{2}$ of a full module $\cM$ are called
{\it associates} if there exists $\eta \in \cU_{\cM}^{\bbK}$ such that
$\mu_{1} = \eta \mu_{2}$.
\end{definition} 

Note that if $\cM$ is a full module, then $\cM^{\bbK}=\cM$.

It is known that there are only finitely many pairwise non-associate elements
with given norm in a full module $\cM$ (see \cite[Corollary to Theorem~5, pg. 90]{BS1}).

While some upper bounds for the number of such non-associates are known
(see the result from \cite{Vout} cited below), it is not known how well these
bounds reflect the actual behaviour of these numbers.

It would not be unreasonable to suspect that this number depends on the field
$\bbK$. In fact, given a result of Schmidt on norm-form equations to
be cited below, one might even believe that this number depends only on $r$.

However, we show here that this is not correct. In particular, we have the following result.

\begin{theorem}
\label{thm:countereg}
For any positive integer $N$ and any number field $\bbK$ with $[\bbK:\bbQ] \geq 3$,
there exists a full module $\cM_{N} \subseteq \bbK$ with at least $N$ pairwise
non-associates of norm $1$.
\end{theorem}

\begin{note}
The full modules, $\cM_{N}$, that we construct here are not rare or exotic in
structure. In fact, it will be apparent in Section~\ref{sect:3} that they are
plentiful and simply-defined -- this is even more striking in Note~\ref{note:more-gen}
there.
\end{note}

\begin{note}
Lemma~\ref{lem:quad} shows that such a result is not true if $\bbK$ is a quadratic
extension of $\bbQ$.

Our construction in this paper fails for quadratic fields, as it should from
Lemma \ref{lem:quad}, since at least three generators of the modules
are required: \\
(i) $1 \in \cM_{N}$,\\
(ii) a fixed unit $\epsilon \in \cM_{N}$ and\\
(iii) a third generator, dependent on $N$, must be in $\cM_{N}$.
\end{note}

\subsection{Norm form equations}

\begin{definition} 
A {\it norm form} $F(\bX) = F \left( X_{1}, \ldots, X_{n} \right)$ is a polynomial
in $\bbQ \left[ X_{1}, \ldots, X_{n} \right]$ that can be expressed as
$$
F(\bX) = a \cN_{\bbK/\bbQ} \left( \al_{1}X_{1} + \cdots + \al_{n}X_{n} \right) 
$$
where $a$ is a non-zero rational number, $\al_{1}, \ldots, \al_{n}$ lie in an
algebraic number field $\bbK$ and $\cN_{\bbK/\bbQ}$ denotes the norm from $\bbK$
to $\bbQ$.
\end{definition} 

For $i=1, \ldots, r$, we let $\sigma_{i}$ denote the isomorphic embeddings
of $\bbK$ into $\bbC$ and write $\al^{(i)}=\sigma_{i}(\al)$ for any $\al \in \bbK$. 
With $L(\bX) = \al_{1}X_{1} + \cdots + \al_{n}X_{n}$, as above, and 
$L^{(i)}(\bX) = \al_{1}^{(i)}X_{1} + \cdots + \al_{n}^{(i)}X_{n}$ for $i=1,\ldots, r$.
We can write $F(\bX)$ in the form 
\begin{equation}
\label{eq:altform}
F(\bX) = a L^{(1)}(\bX) \cdots L^{(r)}(\bX).
\end{equation}

\begin{definition} 
We call two modules $\cL$ and $\cM$ proportional if there is a
fixed $\sigma \neq 0$ such that $\cM = \sigma \cL$.

A module $\cM$ (and hence $F(\bX)$) is called {\it degenerate} 
if it contains a submodule $\cM_{0}$ that is proportional to a full
module $\cL$ in some subfield $\bbL$ of $\bbK$, where $\bbL$ is
neither $\bbQ$ nor an imaginary quadratic field.
\end{definition} 

This definition was formulated by Schmidt \cite{Schm1,Schm2} 
in the early 1970's and he showed that there are non-zero 
rational numbers $m$ such that $F(\bX)=m$ has infinitely many 
solutions in $\bbZ^{n}$ if and only if $F$ is degenerate. 
Moreover, he was also able to show that even if $F$ is 
degenerate, then there is a notion of a family of solutions 
such that there are only finitely many families of solutions 
of $F(\bX)=m$. 

\begin{definition} 
Suppose that $\cN_{\bbK/\bbQ}(\al)=m$ has a solution $\al \in \cM^{\bbL}$, 
then every element of $\al \cU_{\cM}^{\bbL}$ is also a solution and 
we call this set of solutions a {\it family} of solutions. 
Similarly, the set of all elements $\bx \in \bbZ^{n}$ such that 
$L(\bx) \in \al \cU_{\cM}^{\bbL}$ is called a {\it family} of solutions. 
\end{definition} 

That this is a natural notion of a family of solutions is 
probably best seen by means of examples and so we invite the 
reader to consult those presented in \cite[Section~3]{Schm2} 
and \cite[Section~VII.3]{Schm3}.  

By the end of the 1980's, Schmidt had proven his quantitative subspace
theorem \cite{Schm4} and used it to establish upper bounds that depend
only on $m$, $n$ and $r$ for the number of solutions of the norm form 
equation $F(\bX)=m$ when $F$ is non-degenerate \cite{Schm5}.

At that time, Schmidt posed to this author the question of what sort of bounds
one could obtain for the number of families of solutions of degenerate norm forms.  

The author \cite[Theorem~V.1]{Vout} obtained a bound depending only 
on $m,n$ and $r$ for the number of full submodules in subfields of
$\bbK$ such that any solution of $F(\bX)=m$ must lie in the union of 
these submodules. Gy\H{o}ry \cite[Theorem~7]{Gyory} has independently
established this same result. The most recent results in this area are
due to Evertse and Gy\H{o}ry \cite{EG}. Their results are much more
general than the following, but their Theorem~1 implies that the solutions
of $F(\bX)=1$ lie in the union of at most
$$
\left( 2^{33} r^{2} \right)^{n(n+1)(2n+1)/3-2}
$$
full submodules of subfields of $\bbK$.

As was mentioned in the previous paragraph, Gy\H{o}ry's work in this area has
been much more general. He has generalised the concept of a family of solutions
for the norm form setting to that of decomposable form equations.

Decomposable forms include not only norm forms, but also discriminant forms, index forms,
resultant forms, reducible binary forms and other kinds of forms as well. Moreover, Gy\H{o}ry
considers these decomposable form equations over number fields and, more generally, over
finitely generated fields.

Similar to the author's results to be cited in the next paragraph,
in \cite{EG, Gyory}, Gy\H{o}ry obtained explicit upper bounds for the number of
families which depend on certain ``indices'' associated with the module.

For full modules, the author could only establish bounds for the number of
families of solutions which depended more closely on the given module.
In particular, Lemma~V.2 of \cite{Vout} states that if $\cO_{\cM}$ is the
ring of coefficients of  $\cM$ and $[\cU_{\bbK}: \cU_{\cM}]$ is the index of
the unit group of $\cO_{\cM}$ in the unit group of $\bbK$, then the solutions
of the norm form equation $F(\bX) = m$ lie in the union of at most 
\begin{equation}
\label{eq:unit-based-bound}
[\cU_{\bbK}: \cU_{\cM}] { \tau \left( \left| m \right| \right) }^{r}
\end{equation}
families, where $\tau$ is the function which counts the number of positive
divisors of a rational integer. In Theorem~V.2 of \cite{Vout}, an upper bound
in terms of the coefficients of $\cM$ was obtained.

For a full module, $\cM$, the number of families of solutions of the associated norm-form
equation $F(\bX)=m$ is equal to the number of non-associates in $\cM$ of norm
$m/a$. Hence, from Theorem~\ref{thm:countereg}, we obtain the following Corollary.

\begin{corollary}
\label{cor:countereg}
For any positive integer $N$ and any number field $\bbK$ with $[\bbK:\bbQ] \geq 3$,
there exists a full module $\cM_{N} \subseteq \bbK$ such that the equation
$\cN_{\bbK/\bbQ}(\mu) = 1$ has at least $N$ families of solutions with $\mu \in \cM_{N}$.
\end{corollary}

Despite Corollary~\ref{cor:countereg}, for norm forms in three variables which do not
arise from full modules, we are able to get a bound of the desired form. In fact, this is
only possible because Corollary~\ref{cor:countereg} is not true when $\bbK$ is a quadratic field.

\begin{theorem}
\label{thm:family-bound}
Let $\al_{1}, \al_{2}$ and $\al_{3}$ be algebraic numbers 
which are linearly independent over $\bbQ$. Putting  
$\bbK = \bbQ \left( \al_{1}, \al_{2}, \al_{3} \right)$, $r=[\bbK : \bbQ]$ 
and $L(\bX) = \al_{1}X_{1} + \al_{2}X_{2} + \al_{3}X_{3}$,  
we consider the norm form equation 
\begin{equation}
\label{eq:nf}
F(\bX) = a \cN_{\bbK/\bbQ}(L(\bX)) = 1, 
\end{equation} 
where $a$ is a non-zero rational number and $F(\bX) \in \bbZ[\bX]$. 

If $\left[ \bbQ \left( \al_{2}/\al_{1}, \al_{3}/\al_{1} \right): 
\bbQ \right] > 3$, then the solutions of this equation lie in at 
most $10^{969} r^{10}$ families.
\end{theorem}

\begin{note}
The restriction that the $\al_{i}$'s be 
linearly independent is no real restriction, for otherwise 
the norm form equation $F(\bX)=1$ becomes a Thue equation. 
\end{note}

\begin{note}
The condition that 
$\left[ \bbQ \left( \al_{2}/\al_{1}, \al_{3}/\al_{1} \right): 
\bbQ \right] > 3$ ensures that the module generated by $L(\bX)$ 
is not proportional to a full module in any subfield of $\bbK$. 
\end{note}

\begin{note}
With the exception of refinements, this work establishes this aspect of the
behaviour of the number of families of solutions of norm form equations in
three variables.
\end{note}

The method of proof of Theorem~\ref{thm:family-bound} fails when $F(\bX)$ is a norm
form in four (or more) variables satisfying analogous conditions. By Schmidt's Subspace
Theorem, all the solutions of $F(\bX)=1$ correspond to elements of certain
three-dimensional subspaces of $\bbQ^{4}$. If one of these subspaces gives
rise to a full module of rank~$3$, then, from Corollary~\ref{cor:countereg},
it is possible for there to be an arbitrarily large number of families of solutions. 

Several questions and directions for further investigation come to mind.

It would be of considerable diophantine interest to determine the nature of
the dependence of the number of families of solutions on $F(\bX)$ or $\cM$.

From the proof of Theorem~\ref{thm:countereg}, it is clear that some dependence
on $[\cU_{\bbK}: \cU_{\cM}]$ as in \eqref{eq:unit-based-bound} is necessary.

Under what circumstances is the number of families independent of the module $\cM$?

Bombieri and Schmidt \cite{BS2} have shown that $O(r)$ is the correct order of growth
for the number of solutions of Thue equations. What is the correct order of growth in
Theorem~\ref{thm:family-bound}?

\section{Preliminary Lemmas to the Proof of Theorem~$\ref{thm:countereg}$} 

The following is Exercise~5 on page~93 of \cite{BS1}. We include a proof
for completeness.

\begin{lemma}
\label{lem:module-facts}
Let $\cM_{1}$ and $\cM_{2}$ be two full $\bbZ$-modules in $\bbK$.
Then $\cM_{1} \bigcap \cM_{2}$ is a full $\bbZ$-module.
\end{lemma}

\begin{proof}
Let $\{ \be_{1,1}, \ldots, \be_{1,r} \}$ and $\{ \be_{2,1}, \ldots, \be_{2,r} \}$
be sets of generators for $\cM_{1}$ and $\cM_{2}$, respectively. Since these are
full modules, each of these sets of generators forms a basis for $\bbK$ as a
$\bbQ$-vector space. Hence each $\be_{2,i}$ can be expressed as a linear combination
over $\bbQ$ of the $\be_{1,j}$'s. Therefore, there exist least positive integers
$d_{2,i}$, such that $d_{2,i} \be_{2,i} \in \cM_{1}$. Therefore,
$d_{2,i}\be_{2,i} \in \cM_{1} \bigcap \cM_{2}$ for each $i$.

Furthermore, the $d_{2,i}\be_{2,i}$'s are linearly independent over $\bbQ$.
Hence $\cM_{1} \bigcap \cM_{2}$ has $r$ generators and is a full module.
\end{proof}

The next lemma is the key result in establishing Theorem~\ref{thm:countereg}.

\begin{lemma}
\label{lem:more-units}
Let $\epsilon$ be a unit in $\bbK$ with norm $1$ and not a root of unity. For
each positive integer, $i$, let $\cM^{(i)}$ be a full module in $\bbK$ with
$\cO^{(i)}$ as its ring of coefficients and $\cU^{(i)}$ as the units of norm
$1$ in $\cO^{(i)}$. Further, let $\ell_{i}$ be the number of distinct multiplicative
cosets of the form $\epsilon^{v}\cU^{(i)}$.

Suppose that:\\
{\rm (a)} $\cO^{(i)}$ is a proper subset of $\cM^{(i)}$,

\noindent
{\rm (b)} $\epsilon$ is in $\cM^{(i)}$ but not in $\cO^{(i)}$ for each $i$,

\noindent
{\rm (c)} the $\ell_{i}$'s are finite and pairwise relatively prime and all
greater than $1$,

\noindent
{\rm (d)} $\displaystyle \cO_{\bigcap_{i=1}^{N} \cM^{(i)}}=\bigcap_{i=1}^{N} \cO^{(i)}$,
for all positive integers $N$.

Then, for all positive integers $N$, ${\displaystyle \bigcap_{i=1}^{N}} \cM^{(i)}$
is a full module containing at least $2^{N}$ units which are pairwise non-associates.
\end{lemma}

\begin{proof}
From Lemma~\ref{lem:module-facts}, it follows that ${\displaystyle \bigcap_{i=1}^{N}} \cM^{(i)}$
is a full module. So it remains only to prove the statement about the units.

Let $\cS \subseteq \{ 1, \ldots, N \}$.

We define $a_{\cS}$ by
$$
a_{\cS} \equiv \left\{ 
\begin{array}{l l}
  0 \bmod \ell_{i} & \quad \mbox{if $i \in \cS$}\\
  1 \bmod \ell_{i} & \quad \mbox{if $i \not\in \cS$}\\
\end{array} \right.
$$
for each $1 \leq i \leq N$.

By condition~(b), $\epsilon \not\in \cO^{(i)}$, so it follows that $\ell_{i}>1$.

By condition~(c), the $\ell_{i}$'s are relatively prime, so
we can find such an $a_{\cS}$ from the Chinese Remainder Theorem.

Next note that if $i \in \cS$, then $\epsilon^{a_{\cS}} \in \cU^{(i)} \subseteq \cO^{(i)} \subseteq \cM^{(i)}$,
by condition~(a).

If $i \not\in \cS$, then $\epsilon^{a_{\cS}-1} \in \cU^{(i)} \subseteq \cO^{(i)}$.
Since $\epsilon \in \cM^{(i)}$, once again $\epsilon^{a_{\cS}} = \epsilon^{a_{\cS}-1} \cdot \epsilon \in \cM^{(i)}$.

Hence $\epsilon^{a_{\cS}} \in {\displaystyle \bigcap_{i=1}^{N}} \cM^{(i)}$ for each $\cS$.

However, if $\cS \neq \cS'$ are two distinct subsets of $\{ 1, \ldots, N \}$,
then, without loss of generality, there is an $i \in \cS$ such that
$i \not\in \cS'$. Therefore, $\epsilon^{a_{\cS}}/\epsilon^{a_{\cS'}} \not\in \cU^{(i)}$
and hence $\epsilon^{a_{\cS}}/\epsilon^{a_{\cS'}} \not\in {\displaystyle \bigcap_{i=1}^{N}} \cU^{(i)}$.

Since there are $2^{N}$ distinct subsets, $\cS$, there are at least $2^{N}$ such units.

Furthermore, by condition~(d), $\displaystyle \cO_{\bigcap_{i=1}^{N} \cM^{(i)}}=\bigcap_{i=1}^{N} \cO^{(i)}$,
so $\displaystyle \cU_{\bigcap_{i=1}^{N} \cM^{(i)}}=\bigcap_{i=1}^{N} \cU^{(i)}$
and so these units are non-associates in
${\displaystyle \bigcap_{i=1}^{N}} \cM^{(i)}$.
\end{proof}

Now we provide some results about the sorts of full modules that we will use to
construct our examples. We start with our definition and notation for them.

\begin{definition}
Suppose that $\al_{1}$ is an algebraic integer of degree $r_{1}$ over $\bbQ$,
that $\al_{2}$ is of degree $r_{2}$ over $\bbQ \left( \al_{1} \right)$ and that
the minimal polynomial of $\al_{2}$ over $\bbZ \left[ \al_{1} \right]$ is monic.
We let $\bbK=\bbQ \left( \al_{1}, \al_{2} \right)$.

For any positive integer $n$, let $\cM_{n}\left( \al_{1}, \al_{2} \right)$
be the $\bbZ$-module in $\bbK$ generated by $\left\{ \al_{1}^{i} \al_{2}^{j}:
\mbox{ $0 \leq i \leq r_{1}-1$, $0 \leq j \leq r_{2}-1$ with
$(i,j) \neq \left( r_{1}-1, r_{2}-1 \right)$} \right\}$
and $n\al_{1}^{r_{1}-1}\al_{2}^{r_{2}-1}$.

$\cM_{n}\left( \al_{1}, \al_{2} \right)$ is a full module in $\bbK$, so
$\cM_{n}\left( \al_{1}, \al_{2} \right)^{\bbK}=\cM_{n}\left( \al_{1}, \al_{2} \right)$
and we can unambiguously denote $\cO_{\cM_{n}\left( \al_{1}, \al_{2} \right)}^{\bbK}$
by $\cO_{n}\left( \al_{1}, \al_{2} \right)$.
%
\end{definition}

\begin{lemma}
\label{lem:coeff-ring}
{\rm (i)} $\cO_{n}\left( \al_{1}, \al_{2} \right)$ is the order generated as
a $\bbZ$-module by $1$ and
$\left\{ n\al_{1}^{i} \al_{2}^{j} \right\}_{0 \leq i \leq r_{1}-1, 0 \leq j \leq r_{2}-1}$
where $i$ and $j$ are not both $0$.

{\rm (ii)} Let $k_{1}, \ldots, k_{N}$ be positive integers with $K_{N}$ as their
least common multiple. Then
$$
\bigcap_{i=1}^{N} \cM_{k_{i}}\left( \al_{1}, \al_{2} \right)
= \cM_{K_{N}}\left( \al_{1}, \al_{2} \right)
$$
and
$$
\bigcap_{i=1}^{N} \cO_{k_{i}}\left( \al_{1}, \al_{2} \right)
= \cO_{K_{N}}\left( \al_{1}, \al_{2} \right).
$$
\end{lemma}

\begin{proof}
(i) First observe that $1 \in \cO_{n}\left( \al_{1}, \al_{2} \right)$.

Since $\cO_{n}\left( \al_{1}, \al_{2} \right) \subseteq \bbK$, we can
write any element of $\cO_{n}\left( \al_{1}, \al_{2} \right)$ as
$\sum_{i,j} b_{i,j}\al_{1}^{i}\al_{2}^{j}$ with $b_{i,j} \in \bbQ$.

Suppose
$\sum_{i,j} b_{i,j}\al_{1}^{i}\al_{2}^{j} \in \cO_{n}\left( \al_{1}, \al_{2} \right)$
and arrange the terms so that the pairs $(i,j)$ are ordered lexicographically
(i.e., $\left( i_{1}, j_{1} \right)$ is before $\left( i_{2}, j_{2} \right)$ if
$i_{1}<i_{2}$ or if $i_{1}=i_{2}$ and $j_{1}<j_{2}$). Let $\left( i_{0}, j_{0} \right)$
be the last pair such that $b_{i,j} \not\equiv 0 \bmod n$.

If $\left( i_{0}, j_{0} \right)=(0,0)$, then since $b_{0,0} \cdot 1 \in \cM_{n}\left( \al_{1}, \al_{2} \right)$,
we must have $b_{0,0} \in \bbZ$.

If $\left( i_{0}, j_{0} \right) \neq (0,0)$, then
$\left( \sum_{i,j} b_{i,j}\al_{1}^{i}\al_{2}^{j} \right) \cdot
\left( \al_{1}^{r_{1}-1-i_{0}}\al_{2}^{r_{2}-1-j_{0}} \right)
=b_{i_{0},j_{0}}\al_{1}^{r_{1}-1}\al_{2}^{r_{2}-1}$ plus an element of the form $n\cM_{1}\left( \alpha_{1}, \alpha_{2} \right)$
(i.e., in $\cM_{n}\left( \alpha_{1}, \alpha_{2} \right)$) plus ``smaller''
terms lexicographically. This product is also in $\cM_{n}\left( \al_{1}, \al_{2} \right)$
and since $\left\{ \al_{1}^{i}\al_{2}^{j} \right\}_{0 \leq i \leq r_{1}-1, 0 \leq j \leq r_{2}-1}$
form a basis of $\bbK$ so this product has a unique representation of terms of
this basis, we see that $n|b_{i_{0},j_{0}}$.

Proceeding inductively, it follows that $\cO_{n}\left( \alpha_{1}, \alpha_{2} \right)$
is contained in the order generated by $1$ and
$\left\{ n\al_{1}^{i} \al_{2}^{j} \right\}_{0 \leq i \leq r_{1}-1, 0 \leq j \leq r_{2}-1}$
where $i$ and $j$ are not both $0$.

Since the minimal polynomial of $\al_{2}$ over $\bbZ \left[ \al_{1} \right]$ is
monic, it is immediate that
$n\al_{1}^{i}\al_{2}^{j} \cM_{n}\left( \al_{1}, \al_{2} \right)
\subseteq n\cM_{1}\left( \al_{1}, \al_{2} \right) \subseteq \cM_{n}\left( \al_{1}, \al_{2} \right)$
and so (i) follows.

(ii) We prove here a more general result from which both statements follow.

Let $S$ be any subset of ordered pairs of the form $(i,j)$ with $0 \leq i \leq r_{1}-1$
and $0 \leq j \leq r_{2}-1$. Define $\cM_{n,S}\left( \al_{1}, \al_{2} \right)$
to be the full module generated by
$\left\{ \al_{1}^{i}\al_{2}^{j} \right\}_{(i,j) \not\in S}$
and $\left\{ n\al_{1}^{i}\al_{2}^{j} \right\}_{(i,j) \in S}$.

We prove that
$$
\bigcap_{i=1}^{N} \cM_{k_{i},S}\left( \al_{1}, \al_{2} \right)
= \cM_{K_{N},S}\left( \al_{1}, \al_{2} \right).
$$

Note that we need only prove this for $N=2$ as it follows in general, by induction on $N$.

Suppose that $\be \in \cM_{k_{1},S}\left( \al_{1}, \al_{2} \right)
\bigcap \cM_{k_{2},S}\left( \al_{1}, \al_{2} \right)$.
Then
$$
\be = \sum_{(i,j) \not\in S} a_{i,j} \al_{1}^{i}\al_{2}^{j}
+ \sum_{(i,j) \in S} a_{i,j} k_{1}\al_{1}^{i}\al_{2}^{j}
= \sum_{(i,j) \not\in S} b_{i,j} \al_{1}^{i}\al_{2}^{j}
+ \sum_{(i,j) \in S} b_{i,j} k_{2}\al_{1}^{i}\al_{2}^{j},
$$
where the $a_{i,j}$'s and the $b_{i,j}$'s are integers.

Since $\be$ has a unique representation as a linear combination of the
$\al_{1}^{i}\al_{2}^{j}$'s with rational coefficients, it must be the case that
$a_{i,j}=b_{i,j}$ for $(i,j) \not\in S$ and that
$a_{i,j}k_{1}=b_{i,j}k_{2}$ for $(i,j) \in S$. Hence, $k_{2}$ is a divisor of
$a_{i,j}k_{1}$ for $(i,j) \in S$, that is
$a_{i,j}k_{1}=a_{i,j}'\lcm \left( k_{1},k_{2} \right)$ for $(i,j) \in S$.

Thus $\be \in \cM_{K_{2},S}\left( \al_{1}, \al_{2} \right)$ and so
$\cM_{k_{1},S}\left( \al_{1}, \al_{2} \right) \bigcap \cM_{k_{2},S}\left( \al_{1}, \al_{2} \right)
\subseteq \cM_{K_{2},S}\left( \al_{1}, \al_{2} \right)$.

Now we prove the other inclusion. Suppose that $\be \in \cM_{K_{2},S}\left( \al_{1}, \al_{2} \right)$.
Then
$$
\be =
\sum_{(i,j) \not\in S} a_{i,j} \al_{1}^{i}\al_{2}^{j}
+ \sum_{(i,j) \in S} a_{i,j} K_{2}\al_{1}^{i}\al_{2}^{j},
$$
where the $a_{i,j}$'s are integers and hence $k_{1}, k_{2}$ both divide all of
$K_{2}a_{i,j}$ for $(i,j) \in S$. Therefore,
$\be \in \cM_{k_{1},S}\left( \al_{1}, \al_{2} \right) \bigcap \cM_{k_{2},S}\left( \al_{1}, \al_{2} \right)$,
so
$$
\cM_{K_{2},S}\left( \al_{1}, \al_{2} \right)
\subseteq \cM_{k_{1},S}\left( \al_{1}, \al_{2} \right) \bigcap \cM_{k_{2},S}\left( \al_{1}, \al_{2} \right).
$$

Together, these set inclusions show that
$$
\cM_{k_{1},S}\left( \al_{1}, \al_{2} \right) \bigcap \cM_{k_{2},S}\left( \al_{1}, \al_{2} \right)
= \cM_{K_{2},S}\left( \al_{1}, \al_{2} \right).
$$

The result for the $\cM_{k_{i}}\left( \al_{1}, \al_{2} \right)$'s holds by
putting $S=\left( r_{1}-1, r_{2}-1 \right)$, while the result for the
$\cO_{k_{i}}\left( \al_{1}, \al_{2} \right)$'s holds by putting $S=\left\{ (i,j): 0 \leq i \leq r_{1}-1,
0 \leq j \leq r_{2}-1 \right\}$ $-(0,0)$.
\end{proof}

\begin{lemma}
\label{lem:seq}
Let $r$ be a positive integer and put $g_{r}(X)
= \displaystyle\prod_{i=1}^{r} \left( X^{i}-1 \right)$.

There exist positive integers $m_{r}$, $n_{r}$ and $s_{r}$ with
$\gcd \left( m_{r}, s_{r} \right)=1$ such that the following hold.

\noindent
{\rm (i)} $g_{r}(x)/n_{r} \in \bbZ$ for all $x \in \bbZ$ with $x \equiv s_{r} \bmod m_{r}$.

\noindent
{\rm (ii)} There exists an infinite sequence of primes $\left\{ p_{r,i} \right\}$
satisfying $p_{r,i} \equiv s_{r} \bmod m_{r}$ for all $i$ and such that the numbers
$g_{r}\left( p_{r,i} \right)/n_{r}$ $(i=1,2\ldots)$, are pairwise relatively
prime integers each greater than $1$.
\end{lemma}

\begin{proof}
(i) Let $p$ be the least prime number greater than $r+1$, set $n_{r}=g_{r}(p)$
and $m_{r}=n_{r}(r+1)!$. Note that $n_{r}=g_{r}(p) \neq 0$, so we can divide by $n_{r}$
in what follows.

For all $x \equiv p \bmod m_{r}$, $g_{r}(x) \equiv g_{r}(p) \equiv n_{r} \bmod m_{r}$.
Since $m_{r}$ is a multiple of $n_{r}$, $g_{r}(x)/n_{r}$ is an integer for such $x$.
Hence we let $s_{r}=p$.

(ii) Notice that $p$ in the proof of part~(i) does not divide $n_{r}$, since $p$
does not divide $p^{j}-1$ for any $j>0$. In addition, $p>r+1$, so $p$ cannot
divide $(r+1)!$. Together, these two statements imply that $p$ cannot divide
$m_{r}$. Therefore, there are infinitely many primes congruent to
$p \bmod m_{r}$, that is $s_{r} \bmod m_{r}$. Among all such primes, we must show that
there are infinitely many such that $g_{r} \left( p_{r,i} \right)/n_{r}$ are
pairwise relatively prime. We define such a collection of primes inductively.

First, let $p_{r,1}$ be the smallest prime congruent to $s_{r} \bmod m_{r}$ such
that all the real roots of $g_{r}(X)/n_{r}-1$ are less than $p_{r,1}$.

Next suppose that we have a set of primes $p_{r,1}, \ldots, p_{r,N}$ that
satisfy the conditions in the lemma. 

We now find a prime $p_{r,N+1}$ such that $p_{r,1}, \ldots, p_{r,N+1}$ satisfy
the conditions in the lemma.

Let $\cP_{N}$ be the set of all primes that divide 
$$
\Pi_{N} = \prod_{i=1}^{N} \frac{g_{r} \left( p_{r,i} \right)}{n_{r}}.
$$

Let $q \in \cP_{N}$.

Since $g_{r} \left( p_{r,i} \right) \equiv g_{r} \left( s_{r} \right) \equiv g_{r}(p) \equiv n_{r} \bmod m_{r}$
for $p_{r,i} \equiv s_{r} \bmod m_{r}$, $g_{r}\left( p_{r,i} \right)/n_{r} \equiv 1 \bmod (r+1)!$.
So $\gcd \left( g_{r} \left( p_{r,i} \right) / n_{r}, (r+1)! \right) = 1$. Hence
$q > r+1$ or, more conveniently for what follows, $q-1 \geq r+1$.

Now the zeroes of $g_{r}(X) \bmod q$ are the roots of unity $\bmod \, q$ of
order at most $r$. Since $q$ is prime, there always exists a primitive root,
$b_{q}$, modulo $q$, i.e., a number $b_{q}$ such that $b_{q}^{q-1} \equiv 1 \bmod q$
and $b_{q}^{k} \not\equiv 1 \bmod q$ for $0 < k < q-1$. Therefore, $b_{q}$ is a
primitive $q-1$-st root of unity $\bmod q$. Now since
$q-1 \geq r+1$, $g_{r} \left( b_{q} \right) \not\equiv 0 \bmod q$.
Therefore, for each $q \in \cP_{N}$, there is a non-zero congruence class, $b_{q}$,
such that $g_{r} \left( b_{q} \right) \not\equiv 0 \bmod q$.

We choose $p_{r,N+1}$ to be a prime satisfying $p_{r,N+1}>p_{r,N}$, $p_{r,N+1} \equiv s_{r} \bmod m_{r}$
and $p_{r,N+1} \equiv b_{q} \bmod q$ for each $q \in \cP_{N}$. From this last
condition, we have $g_{r} \left( p_{r,N+1} \right) \not\equiv 0 \bmod q$ for any
$q \in \cP_{N}$. By Dirichlet's theorem on primes in arithmetic progressions,
there does exist such a $p_{r,N+1}$ (in fact, there are infinitely many such primes). 

Finally, suppose that $p'$ is a prime which divides $g_{r} \left( p_{r,N+1} \right) / n_{r}$.
Then $g_{r} \left( p_{r,N+1} \right) \equiv 0 \bmod p'$ and thus $p' \not\in \cP_{N}$.
This shows that $g_{r}\left( p_{r,N+1} \right) / n_{r}$ and $\Pi_{N}$ are relatively
prime as desired.

Furthermore, since $p_{r,N}>p_{r,1}$ for all $N \geq 2$ and $p_{r,1}$ is larger
than all the real roots of $g_{r}(X)/n_{r}-1$, our condition that
$g_{r} \left( p_{r,N} \right)/n_{r}>1$ also holds.
\end{proof}

\begin{lemma}
\label{lem:order}
Let $\bbK$ a number field with $r=[\bbK:\bbQ] \geq 2$ containing an order $\cO$.
For any positive integer $n$, let $\cO_{n}$ be the order generated as a
$\bbZ$-module by $1$ and $n\cO$. Let $\eta \in \cO$ be a unit of norm $1$ and
not a root of unity. Put $\epsilon=\eta^{n_{r}}$, using the notation of
Lemma~$\ref{lem:seq}$. For any prime $p$ satisfying $p \equiv s_{r} \bmod m_{r}$
which does not divide $\disc(\cO)$, let $t$ be the least positive integer such
that $\epsilon^{t} \in \cO_{p}$. Then $t$ divides $g_{r}(p)/n_{r}$.
\end{lemma}

\begin{proof}
First recall from Lemma~\ref{lem:seq}(i) that $g_{r}(p)/n_{r}$ is an integer.

The discriminant of $\bbK$ is a divisor of $\disc(\cO)$ and, by assumption, $p$
is not a divisor of $\disc(\cO)$. Therefore $p$ does not ramify in $\bbK$ and
we have
\begin{equation}
\label{eq:p-fact}
\cO_{\bbK}/(p) \cong \cO_{\bbK}/P_{1} \times \cdots \times \cO_{\bbK}/P_{s}
\end{equation}
via the map that takes $x+(p)$ to $(x+P_{1}, \ldots, x+P_{s})$
(see Theorem~2, p.~111 of \cite{Rib}), where the $P_{i}$'s are prime ideals in
$\cO_{\bbK}$, $\cO_{\bbK}/P_{i}$ is a field of cardinality $p^{f_{i}}$ and
$f_{1}+\cdots+f_{s}=r$. The right-hand side of \eqref{eq:p-fact} is a ring
under term-wise addition and multiplication.

Since $\eta$ is a unit, $\eta+P_{i} \neq 0+P_{i}$.
Hence $(\eta+P_{i})^{p^{f_{i}}-1} = 1+P_{i}$ for each $i=1,\ldots,s$.

Let $F=\lcm \left( p^{f_{1}}-1, \ldots, p^{f_{s}}-1 \right)$. Then
$\eta^{F} \equiv 1 \bmod p$. Now since $F|g_{r}(p)$, it follows that
$\eta^{g_{r}(p)} \equiv 1 \bmod p$, which is to say that there exists
$\gamma \in \cO_{\bbK}$ such that  $\eta^{g_{r}(p)}=\epsilon^{g_{r}(p)/n_{r}}=1+p\gamma$.

Letting $\alpha_{1}, \ldots, \alpha_{r}$ be a basis for $\cO$ over $\bbZ$, we
can write
$$
\gamma=\frac{a_{1}\alpha_{1} + a_{2}\alpha_{2} + \cdots + a_{r}\alpha_{r}}{\disc(\cO)},
$$
where $a_{1}, \ldots, a_{r} \in \bbZ$ (see Theorem~9 on page~29 of \cite{Marcus}).

Since $\epsilon^{g_{r}(p)/n_{r}} \in \cO$, we can write
$$
\epsilon^{g_{r}(p)/n_{r}} = 1 + \frac{pa_{1}\alpha_{1}}{\disc(\cO)}
+ \cdots + \frac{pa_{r}\alpha_{r}}{\disc(\cO)}
= b_{1}\alpha_{1} + \cdots + b_{r}\alpha_{r},
$$
where $b_{1}, \ldots, b_{r} \in \bbZ$.

Any such representation must also be unique, since the $\alpha_{i}$'s form a
basis for $\cO$, so we must have $pa_{1}/\disc(\cO), \ldots,
pa_{r}/\disc(\cO) \in \bbZ$. Since, by hypothesis, $p \nmid \disc(\cO)$,
$a_{i}'=a_{i}/\disc(\cO) \in \bbZ$ for each $i$. Therefore,
$\epsilon^{g_{r}(p)/n_{r}} = 1 + pa_{1}'\alpha_{1} + \cdots + pa_{r}'\alpha_{r} \in \cO_{p}$.

The cosets of the form $\epsilon^{v} \cU_{p}$, where $\cU_{p}$ is the group of
units of norm $1$ in $\cO_{p}$, form a group under multiplication. So if
$\epsilon^{t} \in \cU_{p}$, then $t$ must be a divisor of $g_{r}(p)/n_{r}$, as
desired.
\end{proof}

\section{Proof of Theorem~$\ref{thm:countereg}$}
\label{sect:3}

Let $\bbK$ be a number field with $r=[\bbK:\bbQ] \geq 3$ and let $\eta$ be a
unit in $\bbK$ of norm $1$ which is not a root of unity. Put
$\epsilon=\eta^{n_{r}}$, where $n_{r}$ is as in Lemma~\ref{lem:seq}. Let
$\alpha$ be a primitive element of the extension $\bbK/\bbQ(\epsilon)$
whose minimal polynomial over $\bbZ[\epsilon]$ is monic (with $\alpha=1$ if
$\bbK=\bbQ(\epsilon)$).

For any positive integer $N$, let $k_{1}=p_{r,1}, \ldots, k_{N}=p_{r,N}$ be the
first $N$ elements of a sequence of primes satisfying the conditions in Lemmas~\ref{lem:seq}
and \ref{lem:order} (with $\cO=\cM_{1}\left( \epsilon, \alpha \right)$ -- note
that this is an order in $\bbK$ by our conditions on $\alpha$). Put
$K_{N}=k_{1} \cdots k_{N}$.

We are now ready to apply Lemma~\ref{lem:more-units} to prove
Theorem~\ref{thm:countereg}.

We let $\cM^{(i)}=\cM_{k_{i}}(\epsilon,\alpha)$, so that
$\cO^{(i)}=\cO_{k_{i}}(\epsilon,\alpha)$ from Lemma~\ref{lem:coeff-ring}(i).

Notice that $\cO^{(i)}$ is a proper subset of $\cM^{(i)}$, so condition~(a) of
Lemma~\ref{lem:more-units} holds.

Also $\epsilon \in \cM^{(i)}$ and $\epsilon \not\in \cO^{(i)}$, so
condition~(b) of Lemma~\ref{lem:more-units} holds.

Recall from the statement of Lemma~\ref{lem:more-units} that $\ell_{i}$
is the number of distinct multiplicative cosets of the form $\epsilon^{v}\cU^{(i)}$,
where $\cU^{(i)}$ is the group of units of norm $1$ in $\cO^{(i)}$. From
Lemmas~\ref{lem:seq}(ii) and \ref{lem:order}, we know that
$\ell_{i} | \left( g_{r} \left( k_{i} \right)/n_{r} \right)$, which are all pairwise
relatively prime and that $\ell_{i}>1$. Therefore, condition~(c) of
Lemma~\ref{lem:more-units} holds.

Finally, from Lemma~\ref{lem:coeff-ring}(ii), condition~(d) of
Lemma~\ref{lem:more-units} holds.

Since all the conditions in Lemma~\ref{lem:more-units} are satisfied,
$\cM_{K_{N}}(\epsilon,\alpha) = {\displaystyle \bigcap_{i=1}^{N}} \cM_{k_{i}}(\epsilon,\alpha)$
(equality holding by Lemma~\ref{lem:coeff-ring}(ii)) has at least $2^{N}$ units
that are non-associates.

Hence Theorem~\ref{thm:countereg} holds.

\begin{note}
\label{note:more-gen}
These modules, $\cM_{n} \left( \alpha, \epsilon \right)$, are in fact special
cases of more general examples.

Let $\cO$ be any order in $\bbK$,
$\eta$ any unit in $\cO$ of norm $1$ which is not a root of unity and put
$\epsilon=\eta^{n_{r}}$. Let $\varphi: \cO \rightarrow \bbZ$ be any non-trivial
$\bbZ$-module homomorphism such that $\varphi(1)=\varphi(\epsilon)=0$. Define
$\cM_{n, \epsilon, \varphi}$ to be the kernel of the map $\varphi$ $\bmod \, n$
from $\cO$ to $\bbZ/n\bbZ$. It is a full module in $\bbK$.

If $n$ is relatively prime to $\det \left( \varphi \left( \omega_{i} \omega_{j} \right) \right)$,
where $\left\{ \omega_{i} \right\}$ is a basis for $\cO$ as a $\bbZ$-module
$($and the value of this determinant is, in fact, independent of the choice of
basis of $\cO)$, then $\cO_{n, \epsilon, \varphi}$, the ring of coefficients of
$\cM_{n, \epsilon, \varphi}$, is $\bbZ+n\cO$.

Thus the other lemmas in this section can be applied, as here, to construct
modules from these $\cM_{n, \epsilon, \varphi}$'s with arbitrarily many
units that are non-associates.
\end{note}

\section{Preliminary Lemmas to the Proof of Theorem~$\ref{thm:family-bound}$} 

\begin{lemma}
\label{lem:1}
Given $\al_{1}, \ldots, \al_{n}$ which are $\bbQ$-linearly 
independent elements of a number field $\bbK$, let $\cM$ 
be the $\bbZ$-module generated by these $\al_{i}$'s. If 
$n$ is prime and $\bbL$ is a number field such that 
$\cM^{\bbL} = \cM$ then either $\cM$ is proportional to 
a full module in $\bbL$ or $\bbL = \bbQ$.
\end{lemma}

\begin{proof} 
Let $\cM\bbL$ be the set of all products of the form $\al \mu$ 
where $\al \in \bbL$ and $\mu \in \cM$. It is easy to see that 
$\cM\bbL$ is closed under multiplication by elements of $\bbL$. 
Suppose that $\al \mu_{1}, \be \mu_{2} \in \cM \bbL$. Since 
$\cM^{\bbL}=\cM$, there is a non-zero rational integer $a$ such 
that $a \left( \al \mu_{1} + \be \mu_{2} \right) \in \cM$.
Thus $\al \mu_{1} + \be \mu_{2} \in \cM \bbL$. So $\cM \bbL$ 
is also closed under addition and hence is a vector space 
over $\bbL$ of dimension $d$, say. 

Since $\cM^{\bbL} = \cM$, we have $\cM \bbL = \cM\bbQ$ and hence 
$\dim_{\bbQ}(\cM\bbQ) = d[\bbL:\bbQ] = n$. However, $n$ is 
prime so either $d=1$, in which case $\cM$ is proportional 
to a full module in $\bbL$, or $[\bbL:\bbQ]=1$ so that $\bbL = \bbQ$.
\end{proof}

In the case of our theorem, i.e., $n=3$ and $r = [\bbK:\bbQ] > 3$,
Lemma~\ref{lem:1} tells us that if $\cM^{\bbL} = \cM$ then $\bbL = \bbQ$.
This information turns out to be crucial in what follows. 

We will use the heights $H(\cdot)$, $H^{*}(\cdot)$ and $\cH(\cdot)$ 
defined on pages~201 and 204 of \cite{Schm5}. Let
$L=\sum_{j=1}^{n} \al_{j}X_{j}$ be a linear form with coefficients
in an algebraic number field $\bbK$ of degree $r$ over $\bbQ$ and let
$a \in \bbQ^{*}$ be such that the norm form
$$
F(\bX) = a \cN_{\bbK/\bbQ} \left( L (\bX) \right)
= a \prod_{i=1}^{r} \left( \sum_{j=1}^{n} \al_{j}^{(i)} X_{j} \right)
$$
has its coefficients in $\bbZ$. Then the height of F, $H^{*}(F)$, is
defined by
$$
H^{*}(F) = |a| \prod_{i=1}^{r} \left( \sum_{j=1}^{n} \left| \al_{j}^{(i)} \right|^{2} \right)^{1/2}.
$$

According to Lemma~1 of \cite{Schm5}, for the absolute height
$H(L)$ of the linear form $L$, we have
$$
H^{*}(F) = {\rm cont}(F) H(L)^{r},
$$
where ${\rm cont}(F)$ denotes the greatest common divisor of the coefficients
of $F$.

Two norm forms $F$, $G$ are called {\it equivalent}, which we denote
by $F \sim G$, if $G(\bX)=F(B\bX)$ for some matrix $B \in {\rm SL}(n,\bbZ)$.
Now $H^{*}(\cdot)$ is not an invariant under this equivalence, so we define
an invariant height of a norm form $F$, $\cH(F)$, by
$\cH(F) := \min_{G \sim F} H^{*}(G)$, where the minimum is taken over
all norm forms $G$ equivalent to $F$.

To proceed, we now divide the solutions of (\ref{eq:nf}) into large and
small solutions and ``jack up the height'' of the norm form $F(\bX)$.
The point of this last process, which will be explained shortly, is to
replace $F(\bX)$ by a finite number of other norm forms $F_{j}(\bX)$
which are of sufficiently large height so that we can apply known diophantine
techniques to obtain an upper bound on the number of solutions or, when
it works for degenerate norm form equations, families of solutions of 
$F_{j}(\bX)=1$. We create these new forms in such a way, via linear
maps, that the number of solutions (or families of solutions) to the norm
form equation $F(\bX)=1$ is at most the sum of the number of solutions
(or families of solutions) of each of the norm form equations $F_{j}(\bX)=1$.
Since we know the number of such norm forms, we can bound the number
of solutions of $F(\bX)=1$. 

For a prime $p$, let 
$$
A_{0} = \left( \begin{array}{cc}
		       p & 0 \\ 
		       0 & 1 
		 \end{array}
	  \right), 
\hspace{5mm}          
A_{j} = \left( \begin{array}{cc}
		       0 & -1 \\ 
		       p & -j 
		 \end{array}
	  \right)
\mbox{ for $j=1, \ldots, p$.}        
$$

For any $n \geq 2$, we let $E$ be the $(n-2) \times (n-2)$ 
identity matrix and consider the $n \times n$ matrices
$$
B_{j} = \left( \begin{array}{cc}
		     A_{j} & 0 \\ 
		     0     & E 
	       \end{array}
	\right)
\mbox{ for $j=0, \ldots, p$.}        
$$

We can use the linear maps induced by these matrices to express 
$$
\bbZ^{n} = \bigcup_{j=0}^{p} B_{j} \bbZ^{n}. 
$$

For $j=0, \ldots, p$, we put 
$$
F_{j}(\bX) = F \left( B_{j} \bX \right) 
$$
and notice that we can express $F_{j}(\bX)$ in the form 
$$
F_{j}(\bX) = a_{j} \cN_{\bbK/\bbQ} \left( L_{j} (\bX) \right), 
$$
where $a_{j}$ is a non-zero rational number and 
$L_{j}(\bX) = L \left( B_{j} \bX \right)$. 

Moreover, we shall assume that these $F_{j}(\bX)$'s 
are {\it reduced}, that is, $\cH \left( F_{j} \right) 
= H^{*} \left( F_{j} \right)$. This idea comes from 
\cite[p. 208]{Schm5} where it is noted that the 
number of solutions is unaffected by such an assumption.

If we let $p=125000^{3}+21$ (which is prime), then 
\begin{equation}
\label{eq:hbnd}
H \left( L_{j} \right) = \cH { \left( F_{j} \right) }^{1/r} 
		       \geq p^{1/3} > 125000,
\end{equation}
for each $j=0,\ldots,p$, by equations~(5.3) and (5.5) 
of \cite{Schm5}. 

In what follows, we shall drop the subscripts on the $F_{j}$'s 
and $L_{j}$'s in order to simplify our notation. It is also  
at this point where we introduce our definition of small and 
large solutions. 

\begin{definition} 
We define a {\it small solution} of $F(\bx)=1$ to be one with 
\begin{equation}
\label{eq:small}
|\bx| \leq H(L)^{6^{49}r^{3}},
\end{equation}
where $|\bx|$ denotes the ordinary Euclidean absolute value.
A {\it large solution} will be one for which (\ref{eq:small}) does not hold.
\end{definition}

To be able to estimate the number of small solutions 
we need the following lemma. 

\begin{lemma}
\label{lem:lfbnds}
Suppose that $F(\bX)$ is a norm form as in Theorem~$\ref{thm:family-bound}$
which is reduced and satisfies $(\ref{eq:hbnd})$ and that $\bx \in \bbZ^{3}$
is a solution of $F(\bx)=1$. There are three linearly independent forms
$L_{1}(\bX), L_{2}(\bX)$ and $L_{3}(\bX)$ with real coefficients and 
\begin{equation}
\label{eq:lfbnds}
\left| L_{1}(\bx) L_{2}(\bx) L_{3}(\bx) \right| 
< \left| \det \left( L_{1}, L_{2}, L_{3} \right) \right| H(L)^{-2/3},
\end{equation}
where $\det \left( L_{1}, L_{2}, L_{3} \right)$ 
is the determinant of the coefficient matrix. 
\end{lemma}

\begin{proof}
From equation~(6.1) in \cite{Schm5} applied with $n=3$, there exist such
linear forms with 
\begin{eqnarray*}
\left| L_{1}(\bx) L_{2}(\bx) L_{3}(\bx) \right| 
& < & \frac{8^{3}}{(3!)^{1/2}V(3)} \left| \det \left( L_{1}, L_{2}, L_{3} \right) \right| \cH(F)^{-1/r} \\
& < & 50 \left| \det \left( L_{1}, L_{2}, L_{3} \right) \right| \cH(F)^{-1/r},
\end{eqnarray*}
since, $V(3)$, the volume of the unit ball in $\bbR^{3}$, is $4\pi/3$.

Note that, as in Section~6 of \cite{Schm5}, $L$ and its conjugate linear forms
do not necessarily have real coefficients. However, as there, the procedure in
\cite[Section~2]{Schm4} can be applied to obtain the $L_{1}(\bX), L_{2}(\bX)$
and $L_{3}(\bX)$ required here.

Since we have assumed that $F$ is reduced, we know that 
$\cH(F)^{1/r}=H(L)$. By the inequalities in (\ref{eq:hbnd}), 
we know that $50/H(L) < 1/H(L)^{2/3}$ and the lemma follows.
\end{proof}

\begin{lemma}
\label{lem:small}
Suppose that $F(\bX)$ is a reduced norm form as in 
Theorem~$\ref{thm:family-bound}$ which satisfies $(\ref{eq:hbnd})$.  
The small solutions of $F(\bx)=1$ lie in the union of not  more than
$2^{95} 3^{108} r^{9}$ proper linear subspaces of $\bbQ^{3}$.
\end{lemma}

\begin{proof}
Putting $B=H(L)^{6^{49}r^{3}}, P=H(L)^{2/3}$ and
$Q = (\log B)/(\log P) = 2^{48} 3^{50} r^{3}$, we have
$P=H(L)^{2/3} > 2500 > 1296 = (3!)^{4}$, so we can 
apply Schmidt's explicit version of the gap principle 
\cite[Lemma~3.1]{Schm4} to (\ref{eq:lfbnds}) to show 
that for any choice of $L_{1}(\bX), L_{2}(\bX)$ and 
$L_{3}(\bX)$ the solutions of (\ref{eq:lfbnds}) lie at most
$2^{96} 3^{109} r^{6}$ proper linear subspaces of $\bbQ^{3}$.

These $L_{i}(\bX)$'s are obtained from the $L^{(i)}(\bX)$'s and so there are
$\binom{r}{3}$ different ways of choosing them. Consideration of all these
choices leads to the proof of the lemma.
\end{proof}

Let us now turn to the large solutions. 

We first normalise the $L^{(i)}(\bX)$'s: put 
$M_{i}(\bX) = { \left| L^{(i)} \right| }^{-1} L^{(i)}(\bX)$
for $i=1,\ldots,r$. 

\begin{lemma}
\label{lem:lfbndl}
Suppose that $F(\bX)$ is a reduced norm form as in 
Theorem~$\ref{thm:family-bound}$ which satisfies $(\ref{eq:hbnd})$ 
and that $\bx$ is a large solution of $F(\bx)=1$. There are 
integers $1 \leq i_{1} < i_{2} < i_{3} \leq r$ such that 
\begin{equation}
\label{eq:lfbndl}
\left| M_{i_{1}}(\bx) M_{i_{2}}(\bx) M_{i_{3}}(\bx) \right| 
< \left| \det \left( M_{i_{1}}, M_{i_{2}}, M_{i_{3}} \right) \right| 
  |\bx|^{-1/\left( 2 \cdot 6^{48} \right)}. 
\end{equation}
\end{lemma}

\begin{proof}
By Lemma~5 of \cite{E1}, $H(L) \geq \Delta_{\bbK}^{1/(2r(r-1))}$ where
$[\bbK:\bbQ]=r$ and $\Delta_{\bbK}$ denotes the absolute value of the
discriminant of $\bbK$. From Minkowski's theorem (see Corollary~3
on page~137 of \cite{Marcus}), we know that $\Delta_{\bbK} \geq 2$.
Therefore, $H(L)^{10r^{2}} \geq 32>27$ and so, from (\ref{eq:small}),
$|\bx| > H(L)^{12r^{3}+10r^{2}} > 27H(L)^{12r^{3}}$, since $\bx$ is a
large solution of (\ref{eq:nf}).  

Taking $\eta = r (2n)^{-n2^{n+1}}$, which is admissible 
by Lemma~7(ii) of \cite{Schm5}, and $n=3$, it follows from
Lemma~8 of \cite{Schm5} that for every solution $\bx$ of
(\ref{eq:nf}) with $|\bx| > 27H(L)^{12r^{3}}$ there are indices
$1 \leq i_{1} < i_{2} < i_{3} \leq r$ such that
\begin{equation}
\label{eq:lfbndp}
\left| M_{i_{1}}(\bx) M_{i_{2}}(\bx) M_{i_{3}}(\bx) \right| 
< |\bx|^{-1/6^{48}}.
\end{equation}
In fact, Schmidt \cite{Schm5} proved his Lemma~8 under the assumption
that the norm form under consideration is non-degenerate. However,
this assumption was not used in the proof, therefore the lemma, and
hence the inequality, applies as well to degenerate norm forms. 

Let $\bbL$ be the field defined by the property that a subfield 
$\bbF$ of $\bbK$ has $\cM^{\bbF}=\cM$ if and only if $\bbF \subset \bbL$. 
The only place where non-degeneracy is used in the proof of 
Lemma~8 of \cite{Schm5} is on p.214 where it is required that 
$\bbL$ is either $\bbQ$ or an imaginary quadratic field. 

However, we concluded from Lemma~\ref{lem:1} that $\bbL = \bbQ$ 
in our application here and hence the proof of Lemma~8 of 
\cite{Schm5} is also valid here. 

By (5.3) of \cite{Schm4}, 
$$
\left| \det \left( M_{i_{1}}, M_{i_{2}}, M_{i_{3}} \right) \right| 
\geq H(L)^{-3r^{3}}, 
$$
so that, 
\begin{eqnarray*}
\left| M_{i_{1}}(\bx) M_{i_{2}}(\bx) M_{i_{3}}(\bx) \right| 
& < & \left| \det \left( M_{i_{1}}, M_{i_{2}}, M_{i_{3}} \right) \right| 
      H(L)^{3r^{3}} |\bx|^{-1/6^{48}} \\ 
& < & \left| \det \left( M_{i_{1}}, M_{i_{2}}, M_{i_{3}} \right) \right| 
      |\bx|^{-1/ \left( 2 \cdot 6^{48} \right)}, 
\end{eqnarray*}
by our definition of large solutions and (\ref{eq:lfbndp}).
\end{proof}

\begin{note}
It is the lower bound for the determinant in the proof of this lemma, in
particular the exponent on $H(L)$, which dictates our definitions of small
and large solutions. This definition, in turns, affects the number of
subspaces in which the solutions of our norm form equation can belong.
Therefore, an improved lower bound for this determinant would lead to 
an improvement in Theorem~\ref{thm:family-bound}. However, Evertse \cite{E2}  
has shown that in general this bound is best possible with respect to the
exponent on $H(L)$. Hence Theorem~\ref{thm:family-bound} seems to be the
limit of this method in terms of the dependence on $r$. 
\end{note}

To count the families of large solutions, Schmidt used his 
quantitative subspace theorem \cite{Schm4}. Here we shall 
use Evertse's refinement of this result. 

\begin{lemma}
\label{lem:quant}
Let $L_{1}, \ldots, L_{n}$ be linearly independent linear forms 
in $n$ variables such that the field formed by adjoining the 
coefficients of any of the $L_{i}$'s to $\bbQ$ is of degree at 
most $D$ over $\bbQ$ and $H \left( L_{i} \right) \leq H$. 
For every $\delta$ with $0 < \delta < 1$ there are proper linear
subspaces $\cL_{1}, \ldots, \cL_{t}$ of $\bbQ^{n}$ with 
$$
t \leq 2^{60n^{2}} \delta^{-7n} \log (4D) \log \log (4D) 
$$
such that every solution $\bx \in \bbZ^{n}$ of 
\begin{equation}
\label{eq:subineq}
\left| L_{1}(\bx) \cdots L_{n}(\bx) \right| 
< \left| \det \left( L_{1}, \ldots, L_{n} \right) \right| |\bx|^{-\delta} 
\end{equation}
with $\gcd \left( x_{1}, \ldots, x_{n} \right) = 1$ 
and $|\bx| \geq H$ lies in $\cL_{1} \cup \cdots \cup \cL_{t}$. 
\end{lemma}

\begin{proof}
This is the Corollary of \cite{E1}.
\end{proof}

\begin{lemma}
\label{lem:big}
Suppose that $F(\bX)$ is a reduced norm form as in 
Theorem~$\ref{thm:family-bound}$ which satisfies $(\ref{eq:hbnd})$. 
The large solutions of $F(\bX)=1$ lie in at most 
$2^{1570} 3^{1007} r^{3} \log^{2} r$ proper linear subspaces of $\bbQ^{3}$.
\end{lemma}

\begin{proof}
This is a simple consequence of Lemmas~\ref{lem:lfbndl} and \ref{lem:quant},
taking $n=3, \delta=1/(2 \cdot 6^{48}), L_{1}=M_{i_{1}}, L_{2}=M_{i_{2}}$ and
$L_{3}=M_{i_{3}}$. Since these $L_{i}$'s are normalised linear forms coming
from our original $L^{(i)}$'s, we have $D \leq 8r^{3}$ and
$H = \max_{i} H \left( L^{(i)} \right)$, by the product formula. The lemma then
follows from a simple calculation upon noting that $r \geq 4$, that there are
$\binom{r}{3}$ possibilities for $1 \leq i_{1} < i_{2} < i_{3} \leq r$ and that
$r(r-1)(r-2)\log (32r^{3}) \log \log (32r^{3})<4r^{3}\log^{2}r$ for $r \geq 4$.
\end{proof}

\begin{lemma}
\label{lem:quad}
Let $f(X,Y) = aX^{2}+bXY+cY^{2}$ be a norm form with integer 
coefficients then the equations $f(X,Y) = \pm 1$ each have at 
most one family of solutions.
\end{lemma}

\begin{proof}
This follows from Theorem~5 of Section~2.7 of \cite{BS1} and 
the discussion that follows. Theorem~5 states that there is a 
one-to-one correspondence between the families of solutions of 
$f(x,y)=m$ and the modules, $\cA$, in a certain class which 
have norm $m$ and lie in the coefficient ring of the module 
corresponding to $f(x,y)$ (using their definition of classes of 
modules and of the norm of a module). The discussion on p.~144 
demonstrates that finding such modules reduces to the problem
of finding all integers $A$ and $B$ such that $-A \leq B < A$, 
$B^{2}-4AC$ is the discriminant of $f$ and $m=AS^{2}$ for 
some $C,S \in \bbZ$. 

Borevich and Shafarevich show, at the bottom of p.~142 and the 
top of p.~143 of \cite{BS1} that it suffices to consider only 
positive integers $m$. So, putting $m=1$, the last condition 
on $A$ and $B$ shows that $A=1$. Combining this with the first 
condition, we find that either $B=0$ or $B=-1$. This means 
that $D=-4C$ or $D=1-4C$. Clearly only one of these can be true 
and hence there is at most one such module, i.e., at most one 
family of solutions.
\end{proof}

\section{Proof of Theorem~$\ref{thm:family-bound}$} 

Combining Lemmas~\ref{lem:small} and \ref{lem:big}, along 
with our discussion of the relation between solutions of the 
$F_{j}(\bX)=1$ and $F(\bX)=1$, we see that the solutions of 
(\ref{eq:nf}) lie in the union of at most 
$$
\left( 125 \, 000^{3} + 22 \right) 
\left( 2^{95} 3^{108} r^{9} 
       + 2^{1570} 3^{1007} r^{3} \log^{2} r \right) 
< 1.1 \cdot 10^{965} r^{9}
$$
proper linear subspaces of $\bbQ^{3}$, since $r \geq 4$. 

The integer points in any proper linear subspace of $\bbQ^{3}$ can
be parametrised as $\bx = T \by$ where $T$ is a linear map from 
$\bbQ^{2}$ into the subspace which sets up a 1-1 correspondence 
between $\bbZ^{2}$ and the integer points in the subspace. 

Thus, restricting our attention to the integer points of the 
subspaces which arise, our norm form $F(\bX)$ becomes $F(T(\bY))$, 
which is a norm form in two variables with integer coefficients. 
We may also write them as  
$$
F_{1}(\bY) = \cN_{\bbK/bQ} \left( \be_{1} Y_{1} + \be_{2} Y_{2} \right) 
	   = 1. 
$$

We must now consider the fields 
$K_{1} = \bbQ \left( \be_{1}/\be_{2} \right)$. 
If $\left[ \bbK_{1}: \bbQ \right] \geq 3$, then the $\bbZ$-module 
generated by $\be_{1}$ and $\be_{2}$ is not a full module and 
so $F_{1}(\bY)$ is a binary form of degree $r$ which is not a 
power of a binary quadratic form. Thus, as a consequence of 
Theorem~1 of \cite{Stew} (take $\epsilon$ sufficiently large), 
$F_{1}(\bY)=\pm 1$ has at most $5600r$ integer solutions. We 
need to consider $F_{1}(\bY) = \pm 1$, since $F_{1}(\bY)$ might 
be the power of a binary form of lower degree. So it remains to 
consider the case when $\bbK_{1}$ is a quadratic field. But by 
Lemma~\ref{lem:quad}, there are at most two families in this case. 

Thus for each subspace there are at most $5600r$ families of 
solutions. Therefore, by our estimate above for the number of 
proper linear subspaces of $\bbQ^{3}$ into which the solutions
must fall, there are at most $1.1 \cdot 10^{965} \cdot 5600r^{10}
< 10^{969} r^{10}$ families of solutions to (\ref{eq:nf}). 

\section{Acknowledgements}

This work originated during the author's stays at the University of Colorado in
Boulder during the 1990's, first as a graduate student and again as a visiting
professor. The author is extremely grateful to Wolfgang Schmidt for his kindness,
generosity and support during these periods, financially as well as with his
time and ideas.

The author also thanks the referees for their careful reading of this paper and
their helpful suggestions. Their advice improved both the results as well as
the presentation here and is appreciated.

\end{document}